 \documentclass[final,1p,times]{elsarticle}

\usepackage{epsfig}

\usepackage{amssymb}
\usepackage{amsthm}
\usepackage{amsmath}
\usepackage{amsfonts}
\usepackage{mathtools}
\usepackage{booktabs,lipsum}
\usepackage{multirow}
\usepackage[colorlinks=true,linkcolor=red]{hyperref}

\usepackage{lineno}
\biboptions{sort&compress}

\journal{??}

\begin{document}

\begin{frontmatter}



\title{Application of Meshfree Method Based on Compactly Supported Radial Basis Function for Solving Unsteady Isothermal Gas Through a Micro-Nano Porous Medium}


\author[mymainaddress]{K. Parand\corref{mycorrespondingauthor}}
\cortext[mycorrespondingauthor]{Corresponding author. Member of research group of Scientific Computing. Fax:
+98 2122431653.}
\ead{k\_parand@sbu.ac.ir}

\author[mymainaddress]{M. Hemami}
\ead{mohammadhemami@yahoo.com}

\address[mymainaddress]{Department of Computer Sciences, Shahid Beheshti University, G.C. Tehran 19697-64166, Iran}

\begin{abstract}

In this paper, we have applied the Meshless method based compactly supported radial basis function collocation   for obtaining the numerical solution of unsteady gas equation. The unsteady gas equation is a second order non-linear two-point boundary value ordinary differential equation on the semi-infinite domain, with a boundary condition in the infinite. The compactly supported radial basis function collocation method reduces the solution of the equation to the solution of a system of algebraic equation. also, we compare the results of this work with some results. It is found that our results agree well with those by the numerical method, which verifies the validity of the present work.

\end{abstract}


\end{frontmatter}


 \section{Introduction}
\label{intro}
\subsection{Unsteady gas problem}
\label{sec:1}
In the study of the unsteady flow of gas through a semi infinite porous medium \citep{Kidder} initially filled with gas at a uniform pressure $P_0 \geq 0$, at time $t=0$, the pressure at the outflow face is suddenly reduced from $P_0$ to $P_1 \geq 0$ ($P_1 =0$ is the case of diffusion into a vacuum) and is, thereafter, maintained at this lower pressure. The unsteady isothermal flow of gas is described by a non-linear partial differential equation 
\begin{equation}
\nabla ^2(P^2) = 2 A\frac{\partial P}{\partial t},
\end{equation}
where the constant $A$ is given by the properties of the medium.In the one dimensional medium extending from $z=0$ to $z= \infty$, this reduces to 
\begin{equation}
\frac{\partial}{\partial z}(P\frac{\partial P}{\partial z})=A \frac{\partial P}{\partial t},
\end{equation}
with the boundary conditions 
\begin{eqnarray}
P(z,0)=P_0,   0<z<\infty; \\
P(0,t)=P_1(<P_0),   0\leq t <\infty.
\end{eqnarray}
To obtain a similarity solution, Authors\citep{Agrawal} introduced the new independent variable  
\begin{equation}
x=\frac{z}{\sqrt{t}}\sqrt{\frac{A}{4P_0}},
\end{equation}
and the dimension-free dependent variable $y$, defined by 
\begin{equation}
y(x)=\frac{1}{\alpha}(1-\frac{P^2(z)}{P^2_0}),
\end{equation}
where $\alpha  = 1- \frac{P_1^2}{P_0^2}$. In terms of the new variable, the problem takes the form (unsteady gas equation) 
\begin{eqnarray}
\frac{d^2 y}{dx^2}+\frac{2 x}{\sqrt{1-\alpha y(x)}}\frac{dy}{dx}=0,\\
\label{Main_EQ}
x>0,~~~~0\leq \alpha \leq 1.
\end{eqnarray}
The typical boundary conditions imposed by the physical properties are 
\begin{equation}
y(0)=1,~~~~y(\infty)=0.
\end{equation} 
A substantial amount of numerical and analytical work has invested so far \citep{Kidder,Roberts} in this model. The main reason of this interest is that the approximation can be used for many engineering purpose. As stated before, the problem of Eq. (\ref{Main_EQ}) was handled by Kidder\citep{Kidder} where a perturbation technique is carried out to include terms of the second order. Wazwaz \cite{Wazwaz} solved this equation non-linearly by modifying the decomposition method and Pad\'e approximation. Parand \citep{Parand1,Parand2,Parand3}  also applied the Lagrangian method, generalized Laguerre polynomials and Bessel collocation method for solving unsteady gas equation.
Rad and Parand  \citep{Rad} proposed an Analytical solution Unsteady Gas equation by Homotopy perturbation method (HPM). 
Kazem  \citep{Kazem} applied the radial basis function (RBF) collocation method  and  Aslam noor \citep{Noor} applied the variational iteration method (VIM) for solving non-linear this equation. khan \citep{Khan} propose a new approach to solve unsteady gas equation.  he applied the modified Laplace decomposition method (MLDM) coupled with Pad\'e approximation to compute a series solution of unsteady flow of gas through a porous medium. Upadhyay and Rai \citep{Upad} using the Legendre wavelet collocation method (yLWCM) to solve this equation.

\label{sec:4}
\subsection{CSRBF}
\label{ssec:41}
Many problems in science and engineering arise in infinite and semi-infinite domains. Different numerical methods have been proposed for solving problems on various domains such as FEM\citep{DH10,DH11}, FDM\citep{DH9,DH10} and Spectral\citep{DH5,Parand2} methods and meshfree method\citep{DH3,DH4}.
The use of the RBF is an one of the popular meshfree method for solving the differential equations  \citep{DH1,DH2}. For many years the global radial basis functions such as Gaussian, Multi quadric, Thin plate spline, Inverse multiqudric and etc was used  \citep{DH6,DH7,DH8}. These functions are globally supported and generate a system of equations  with ill-condition full matrix.To convert the ill-condition matrix to a well-condition matrix, CSRBFs can be used  instead of global RBFs\citep{Shen}. CSRBFs can convert the global scheme into a local one with banded matrices, Which makes the RBF method more feasible for solving large-scale problem \citep{Wong}.
\subsubsection{Wendland's functions}
\label{ssec:42}
The most popular family of CSRBF are Wendland functions. This function introduced by Holger Wendland in 1995 \citep{Wendland,Wendland2}. Wendland starts with the truncated power function $\phi_l(r)=(1-r)^l_+$ which be strictly positive definite and radial on $\mathbb{R}^s$ for $l\geq \lfloor\frac{s}{2}\rfloor+1$ , and then he walks through dimension by repeatedly applying the operator I.\\
\textbf{Definition \citep{MeshfreeFasshhauer}}~~
with $\phi_l(r)=(1-r)^l_+$ we define
 \begin{equation}
 \phi_{s,k}=I^k\phi_{\lfloor\frac{s}{2}\rfloor +k+1},
 \end{equation}
 it turns out that the functions $\phi_{s,k}$ are all supported on [0,1].\\
\textbf{Theorem 1 \citep{MeshfreeFasshhauer}}~~
 The function $\phi_{s,k}$ are strictly positive definite (SPD)  and radial on $\mathbb{R}^s$ and are of the form
 \begin{equation}
\phi _{s,k}(r)=\begin{cases}
p_{s,k}(r)&   r\in [0,1],\\
0&       r>1,
\end{cases}
\end{equation}
with a univariate polynomial $p_{s,k}$ of degree $\lfloor \frac{s}{2}\rfloor+3k+1$. Moreover, ،$\phi_{s,k}\in C^{2k}(R)$ are unique up to a constant factor, and the polynomial degree is minimal for given space dimension $s$ and smoothness $2k$ \citep{MeshfreeFasshhauer}.
Wendland gave recursive formulas for the functions $\phi_{s,k}$ for all $s, k$. We instead list the explicit formulas of \citep{Fasshhauer}.\\
\textbf{Theorem 2 \citep{MeshfreeFasshhauer}}
~~~~The function $\phi_{s,k}$, $k=0,~1,~ 2,~ 3,$ have form
 \begin{eqnarray}
&&\phi_{s,0}=(1-r)^l_+,\\
&&\phi_{s,1}\doteq(1-r) _+^{l+1}[(l+1)r+1],\\
&&\phi_{s,2}\doteq(1-r)_+^{l+2}[(l^2+4l+3)r^2+(3l+6)r+3],\\
&&\phi_{s,3}\doteq(1-r)_+^{l+3}[(l^3+9l^2+23l+15)r^3+(6l^2+36l+45)r^2\\
&&+(15l+45)r+15],
\end{eqnarray}
where $l=\lfloor \frac{s}{2}\rfloor +k+1$, and the symbol $\doteq$ denotes equality up to a multiplicative positive constant.\\
~~~~~~The case $k=0$ follows directly from the definition. application of the definition for the case $k=1$ yields
\begin{eqnarray}
&&\phi_{s,1}(r)=(I\phi_l)(r)
=\int^\infty_r t\phi_l(t)dt  \\
&&=\int^\infty _r t(1-t)^l_+dt
=\int^1_r t(1-t)^ldt\\
&&=\frac{1}{(l+1)(l+2)} (1-r)^{l+1}[(l+1)r+1],
\end{eqnarray}
where the compact support of $\phi_l$ reduces the improper integral to a definite integral which can be evaluated using integration by parts. The other two cases are obtained similarly by repeated application of $I$.\citep{MeshfreeFasshhauer}
We showed the most of wendland functions in Table \ref{Table. 1.}.
\begin{table}[htbp]
\caption{Wendland's compactly supported radial function for various choices of k and s=3.}
\label{Table. 1.}       
\centering  \begin{tabular}{lll}
\hline\noalign{\smallskip}
 $\phi_{s,k}$ & smoothness & SPD   \\
\noalign{\smallskip}\hline\noalign{\smallskip}
 $\phi_{3,0}(r) =(1-r)^2_+$ & $C^0$ & $\mathbb{R}^3$  \\
 $\phi_{3,1}(r)\doteq(1-r)^4_+(4r+1)$ & $C^2$ & $\mathbb{R}^3$ \\
 $\phi_{3,2}(r)\doteq(1-r)^6_+(35r^2+18r+3)$ & $C^4$ & $\mathbb{R}^3$  \\
 $\phi_{3,3}(r)\doteq(1-r)^8_+(32r^3+25r^2+8r+1)$ & $C^6$ & $\mathbb{R}^3$  \\
 $\phi_{3,4}(r)\doteq(1-r)^{10}_+(429r^4+450r^3+210r^2+50r+5)$ & $C^8$ & $\mathbb{R}^3$ \\
 $\phi_{3,5}(r)\doteq(1-r)^{12}_+(2048r^5+2697r^4+1644r^3+566r^2+108r+9)$ & $C^{10}$ & $\mathbb{R}^3$  \\
\noalign{\smallskip}\hline
\end{tabular}
\end{table}

\subsubsection{Wu's functions}
In \citep{Wu} we find another way to construct strictly positive definite radial function compact support. Wu starts with the function
\begin{equation}
\phi (r)=(1-r^2)_+^l,    l\in N,
\end{equation}
which in itself is not positive definite, however, Wu then uses convolution to construct another function that is strictly positive definite and radial on $\mathbb{R}$
\begin{equation}
\phi_l (r)=\int _{-1}^1 (1-t^2)^2_+(1-(2r-t)^2)^l_+dt.
\end{equation}
then, he constructed functions by "dimension walk" using the $D$ operator.\\
\textbf{Definition 2}
~~With $\phi_l(r)=((1-.^2)_+^l*(1-.^2)_+^l)(2r)$ we define\\
\begin{equation}
\phi_{k,l}=D^k\phi_l.
\end{equation}
The function $\phi_{k,l}$ are strictly positive definite and radial on $\mathbb{R}^s$ for $s\leq  2k+1$. are polynomials of degree $4l-2k+1$ on their support 
and in $C^{2(l-k)}$ in the interior of the support. On the boundary the smoothness increases to $C^{2l-k}$. \citep{MeshfreeFasshhauer}\\
a number of Wu's function showed in Table \ref{Table. 2.}.\\

\begin{table}[htbp]
\caption{Wu's compactly supported radial function for various choice of k and l=3.}
\label{Table. 2.}       
\centering  \begin{tabular}{lll}
\hline\noalign{\smallskip}
 $\phi_{s,k}$ & smoothness & SPD   \\
\noalign{\smallskip}\hline\noalign{\smallskip}
$\phi_{0,3}(r)=(1-r)^7_+(5+35r+101r^2+147r^3+101r^4+35r^5+5r^6$&$C^6$&$\mathbb{R}^3$\\
$\phi_{1,3}(r)\doteq(1-r)^6_+(6+36r+82r^2+72r^3+30r^4+5r^5)$&$C^4$&$\mathbb{R}^3$\\
$\phi_{2,3}(r)\doteq(1-r)^5_+(8+40r+48r^2+35r^3+5r^4)$&$C^2$&$\mathbb{R}^3$\\
$\phi_{3,3}(r)\doteq(1-r)^4_+(16+29r+20r^2+5r^3)$&$C^0$&$\mathbb{R}^3$\\
\noalign{\smallskip}\hline
\end{tabular}
\end{table}

The Oscillator \citep{MeshfreeFasshhauer,Gneiting} and Buhman \citep{Buhman} functions are the other kind of CSRBFs can be showed in Tables \ref{Table. 3.} and \ref{Table. 4.}.

\begin{table}[htbp]
\caption{Oscillator's compactly supported radial function for various choice of k.}
\label{Table. 3.}       
\centering  \begin{tabular}{lll}
\hline\noalign{\smallskip}
 $\phi_{k}$ & smoothness & SPD   \\
\noalign{\smallskip}\hline\noalign{\smallskip}
$\phi_{1}(r)=(1-r)^4_+(1+4r-15r^2)$&$C^2$&$\mathbb{R}^3$\\
$\phi_{2}(r)\doteq(1-r)^6_+(3+18r+3r^2-192r^3)$&$C^4$&$\mathbb{R}^3$\\
$\phi_{3}(r)\doteq(1-r)^8_+(15+120r+210r^2-840r^3-3465r^4)$&$C^6$&$\mathbb{R}^3$\\
\noalign{\smallskip}\hline
\end{tabular}
\end{table}

\begin{table}[htbp]
\caption{Buhman's compactly supported radial function.}
\label{Table. 4.}       
\centering  \begin{tabular}{lll}
\hline\noalign{\smallskip}
 $\phi_{n}$ & smoothness & SPD   \\
\noalign{\smallskip}\hline\noalign{\smallskip}
$\phi_1(r)=12r^4\log r-21r^4+32r^3-12r^2+1$&$C^2$&$\mathbb{R}^3$\\
$\phi_2(r)=2r^4\log r-\dfrac{7}{2}r^4+\dfrac{16}{3}r^3-2r^2+\dfrac{1}{6}$&$C^2$&$\mathbb{R}^3$\\
$\phi_3(r)=r^8-\dfrac{84}{5}r^6+\dfrac{1024}{5}r^{\frac{9}{2}}-378r^4+\dfrac{1024}{5}r^{\frac{7}{2}}-\dfrac{84}{5}r^2+1$&$C^3$&$\mathbb{R}^3$\\
$\phi_4(r)=\dfrac{99}{35}r^8-132r^6+\dfrac{9216}{35}r^{\frac{11}{2}}-\dfrac{11264}{35}r^{\frac{9}{2}}+198r^4-\dfrac{396}{35}r^2+1$&$C^4$&$\mathbb{R}^3$\\
\noalign{\smallskip}\hline
\end{tabular}
\end{table}

\section{ CSRBF method}
\subsection{Interpolation by CSRBFs}
\label{ssec:43}
The one-dimensional function $y(x)$ to be interpolated or approximated can be represented by an CSRBF as
\begin{equation}
y(x)\approx y_n(x)=\sum_{i=1}^N \xi_i \phi_i (x)=\Phi ^T(x)\Xi,
\label{systmat}
\end{equation}
where
\begin{eqnarray}
&&\phi_i (x)=\phi (\dfrac{\|x-x_i\|}{r_\omega}),\\
&&\Phi ^T(x)=[\phi_1 (x),\phi_2 (x),\cdots ,\phi _N(x)],\\
&&\Xi =[\xi_1,\xi_2, \cdots , \xi_N ]^T,
\end{eqnarray}
\begin{equation}
\end{equation}
$x$ is the input, $r_\omega$ is the local support domain and $\xi_i$s are the set of coefficients to be determined. By using the local support domain, we mapped the domain of problem to CSRBF local domain. By choosing $N$ interpolate points $(x_j,~ j=1, 2,\cdots,  N)$ in domain:
\begin{equation}
y_j=\sum_{i=1}^N\xi_i\phi_i (x_j)   (j=1, 2, \cdots, N).
\end{equation}
To summarize the discussion on the coefficients matrix, we define
\begin{equation}
A\Xi = Y,
\label{system}
\end{equation}
where :
\begin{eqnarray}
&&Y=[y_1, y_2, \cdots,  y_N]^T,\\
&&A=[\Phi ^T(x_1), \Phi ^T(x_2), \cdots , \Phi ^T(x_N)]^T\\
&&=\begin{pmatrix} 
\phi_1(x_1)&\phi_2(x_1)&\cdots &\phi_N(x_1)\cr \phi_1(x_2)&\phi_2(x_2)&\cdots &\phi_N(x_2)\cr \vdots &\vdots &\ddots &\vdots \cr \phi_1(x_N)&\phi_2(x_N)&\cdots &\phi_N(x_N)
\end{pmatrix}.
\end{eqnarray}
Note that $\phi_i(x_j)=\phi(\dfrac{\|x_i-x_j\|}{r_\omega})$, by solving the system (\ref{system}), the unknown coefficients $\xi_i$  will be achieved.
\subsection{Solving the  transformed Model}
\begin{figure}[htbp!]
\centering \includegraphics[scale=0.4]{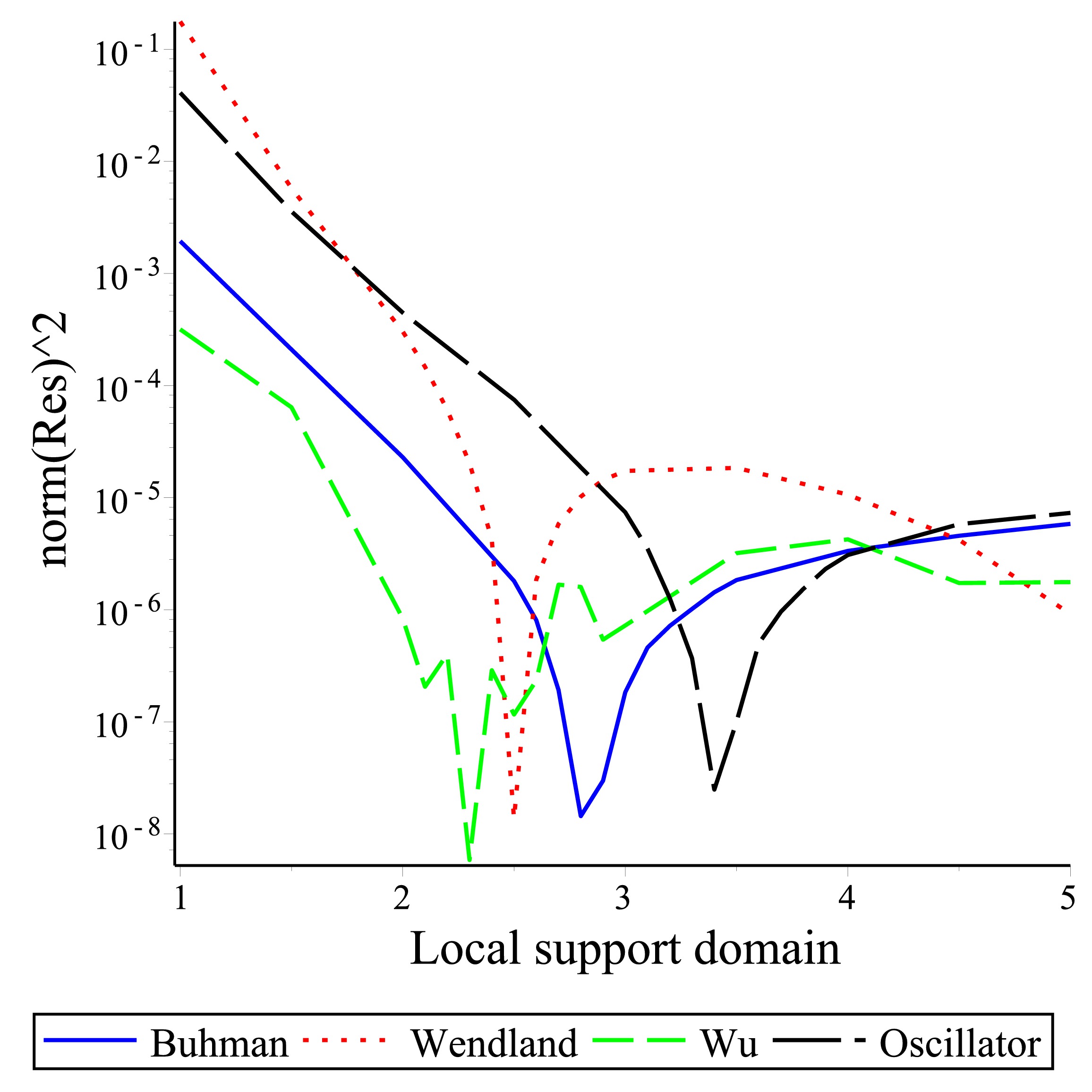}
\caption{The graph of residual error versus local support domain when $N=30$, $\alpha=0.5$ and $xi=1.5$}
\label{figure1}
\end{figure}

In this section, by defining $u(x)$ as a below form and transforming $\frac{dy}{dx}$ and $\frac{d^2y}{dx^2}$ in terms of it,
\begin{eqnarray}
u(x)=1-\alpha y(x) \\
\frac{dy}{dx}= -\frac{du}{dx}\alpha^{-1}\\
 \frac{d^2y}{dx^2}=-\frac{d^2u}{dx^2}\alpha^{-1}
\end{eqnarray}
the problem of unsteady gas in a semi-infinite porous medium (\ref{Main_EQ}) takes the following form :
\begin{equation}
\label{sub4}
\frac{d^2u}{dx^2}+\frac{2x\alpha}{\sqrt{u}}\frac{du}{dx}=0,~~~~x>0, 0\leq \alpha \leq 1,
\end{equation}
and the conditions change to 
\begin{equation}
\label{condition}
u(0)=1-\alpha,~~~~ u(\infty)=1.
\end{equation}
Now we approximate $\frac{du}{dx}$ and $\frac{d^2u}{dx^2}$ as 
\begin{eqnarray}
\label{sub1}
\frac{du}{dx}\simeq \frac{du_n}{dx}=\sum_{i=0}^N \xi_i \phi_i(x),\\
\label{sub2}
\frac{d^2u}{dx^2}\simeq \frac{d^2u_n}{dx^2}=\sum_{i=0}^N \xi_i\phi'(x).
\end{eqnarray} 
By using integral operation $u(x)$ is obtained as
\begin{eqnarray}
\label{sub3}
\int _0^x \frac{du_n}{dr}dr = \sum_{i=0}^N \xi_i  \int _0 ^x \phi_i(r) dr,\\
u_n(x)-u(0)=\sum_{i=0}^N \xi_i  \int _0 ^x \phi_i(r) dr,\\
u(x)\simeq u_n(x)=(\sum_{i=0}^N \xi_i  \int _0 ^x \phi_i(r) dr )+(1-\alpha),
\end{eqnarray}
By substituting (\ref{sub1}), (\ref{sub2}), and (\ref{sub3}) in (\ref{sub4}), we define residual function 
\begin{equation}
Res(x)=\frac{d^2u_n}{dx^2}+\frac{2x}{\sqrt{u_n}}\frac{du_n}{dx}.
\end{equation}
Now, by using $N-1$ interpolate points $\{x_j\}^{N-1}_{i=1}$ plus a condition (\ref{condition}) we can solve the set of equations and consequently, the coefficients $\{\xi_i\}_{i=1}^N$ will be obtained: 
\begin{equation}
\begin{cases}
Res(x_j)=0,~~~~ j=1,2,...,N-1,\\
u_n(L_{\infty})=1.
\end{cases}
\end{equation}.
Collocation point are chosen by :
\begin{equation}
L_{\infty}(\frac{j}{N})^\varrho,~~~j=1,2,...,N-1,
\end{equation}

where $\varrho$ is a arbitrary parameter. Here, we choose $L_\infty=5$ which satisfies $y(L_\infty)<\epsilon$ with $\epsilon$ as a small positive value. The parameter $\varrho$ and local support domain $r_\omega$ that must be selected by the user. But here, by the meaning of residual function, we try to minimize $\|Res(x)\|^2$ by choosing a good local support domain and $\varrho$. We define $\|Res(x)\|^2$ as
\begin{equation}
\|Res(x)\|^2=\int_0^{L_{\infty}} Res^2(x)dx\simeq \sum_{j=0}^m \vartheta_j Res^2(\frac{L_\infty}{2}s_j+\frac{L_\infty}{2}),
\end{equation}
where 
\begin{eqnarray}
\vartheta_j=\frac{L_\infty}{(1-s_j^2)(\frac{d}{ds}P_{m+1}(s)|s=s_j)},~~j=0,1,...,m,\\
P_{m+1}(s_j)=0,~~j=0,1,...,m,
\end{eqnarray}
$P_{m+1}(s)$ is the $(m+1)$th-order Legendre polynomial. The Figure \ref{figure1}. show the minimum of $\|Res(x)\|^2$ which is obtained with $r_\omega$  for cases of$Wendland_{3,5}$, $Wu_{3,3}$, $Buhman_4$ and $Oscillator_3$. 
The result of this section can be summarized in the following algorithm for the BVP:
\begin{equation}
F(x,y(x),\frac{dy}{dx},\frac{d^2y}{dx^2})=0,~~~~y(0)=a,~~~y(\infty)=b.
\end{equation}
\textbf{Algorithm} The algorithm works in the following manner:\\
(1) Choose $N$ center points $\{X_j\}^N_{j=0}$ from domain $[0, L_\infty]$.\\
(2) Approximate $\frac{dy}{dx}$ as the from $\frac{dy}{dx}=\sum_{i=1}^N \xi_i\phi_i(r)$.\\
(3) Obtain $y(x)$ by using defined integral operation $I_\chi(h(x))=\int_0^x h(t)dt$ in the form $u_N(x)= \sum_{i=1}^N \xi_i \int _0^x \phi_i(t) dt+a$\\
(4) Substitiute $u_N(x)$, $\frac{du_N}{dx}$ and $\frac{d^2 u_N}{dx}$ into the main problem and create residual function Res(x).\\
(5) Substitiute collocation points $\{X_j\}_{j=0}^{N-1}$ into the Res(x), along with a boundary condition $y(L_\infty)=b$ and create $N$ equations.\\
(6) Solve the $N$ equations with $N$ unknown coefficients $\{\xi_i\}_{i=1}^N$ and  find the numerical solution.

\section{Results and discussion}
In this section, we compare the applied result of CSRBF with RBF methods \citep{Kazem} and yLWCM \citep{Upad} .\\
In the physical observation of the unsteady gas problem, $y'(0)$ has an important issue \citep{Kidder}.Table \ref{Table. 5.}  and \ref{Table. 6.}.  presents a comparison between the values of $y'(0)$ and $y(x)$ for $\alpha=0.5$ and obtained by the CSRBF and, RBFs and yLWCM.

\begin{table}[htbp] 
\caption{Comparison of  $y(x)$ for $\alpha=0.5$ with $\varrho=1.5$, $N=30$.}
\label{Table. 5.} 
\centering \scriptsize \begin{tabular}{lllllllll}
\hline\noalign{\smallskip}
x & $Wendland_{3,5}$  $r_\omega=2.5$ & $Wu_{3,3}$ $r_\omega=2.3$ & $Oscillator_3$ $r_\omega=3.4$ & $Buhman_4$ $r_\omega=2.8$ & RBF.G & RBF.T& yLWCM \\
\noalign{\smallskip}\hline\noalign{\smallskip}
0.1 & 0.88136427 & 0.88136428 & 0.88136409 & 0.88136468 & 0.88139802 & 0.88137298 & 0.88147552\\
0.2 & 0.76582809 & 0.76582823 & 0.76582774 & 0.76582852 & 0.76588029 & 0.76579834 & 0.76661101\\
0.3 & 0.65599963 & 0.65999935 & 0.65599915 & 0.65600036 & 0.65606928 & 0.65590792 & 0.65727018\\
0.4 & 0.55389758 & 0.55389797 & 0.55389693 & 0.55389814 & 0.55399431 & 0.55375746 & 0.55567752\\
0.5 & 0.46094112 & 0.46094164 & 0.46094037 & 0.46094191 & 0.46107202 & 0.46078383 & 0.46325882\\
0.6 & 0.37797968 & 0.37798027 & 0.37797895 & 0.37798037 & 0.37814380 & 0.37783286 & 0.38078076\\
0.7 & 0.30535020 & 0.30535087 & 0.30534921 & 0.30535091 & 0.30554174 & 0.30522806 & 0.30853930\\
0.8 & 0.24295205 & 0.24295279 & 0.24295099 & 0.24295286 & 0.24316554 & 0.24285666 & 0.24644872\\
0.9 & 0.19033171 & 0.19033246 & 0.19033052 & 0.19033235 & 0.19056419 & 0.19025895 & 0.19409149\\
1.0 & 0.14677064 & 0.14677143 & 0.14677694 & 0.14677139 & 0.14702048 & 0.14671538 & 0.15877555\\
1.2 & 0.08313306 & 0.08313387 & 0.08313169 & 0.08313370 & 0.08341060 & 0.08310227 & --\\
1.4 & 0.04404326 & 0.04404412 & 0.04404186 & 0.04404399 & 0.04433767 & 0.04402800 & --\\
1.8 & 0.01002990 & 0.01003085 & 0.01002839 & 0.01003049 & 0.01034206 & 0.01002778 & --\\
2.2 & 0.00170828 & 0.00170909 & 0.00170673 & 0.00170889 & 0.00202559 & 0.00170941 & --\\
2.6 & 0.00021351 & 0.00021440 & 0.00021195 & 0.00021413 & 0.00053177 & 0.00021524 & --\\
3.0 & 0.00001696 & 0.00001786 & 0.00001540 & 0.00001758 & 0.00033532 & 0.00001878 & --\\
\hline
\hline
$\|Res\|^2$ & 1.3702e-08 & 5.8431e-9 & 2.4820e-08 & 1.4394e-08 & 8.52e-07 & -- & --\\ 
\noalign{\smallskip}\hline
\end{tabular}
\end{table}
\begin{table}[htbp] 
\caption{Comparison of initial slope $y'(0)$ for $\alpha=0.5$ with $\varrho=1.5$, $N=30$.}
\label{Table. 6.} 
\centering \scriptsize \begin{tabular}{lllllllll}
\hline\noalign{\smallskip}
x & $Wendland_{3,5}$  $r_\omega=2.5$ & $Wu_{3,3}$ $r_\omega=2.3$ & $Oscillator_3$ $r_\omega=3.4$ & $Buhman_4$ $r_\omega=2.8$ & RBF.G & RBF.T & yLWCM \\
\noalign{\smallskip}\hline\noalign{\smallskip}
y'(0) & -1.191796 & -1.191806 & -1.191800 & -1.191768 & -1.191498 &  -1.191243 & -1.199258 \\ 
\noalign{\smallskip}\hline
\end{tabular}
\end{table}

 In addition, the value of $y'(0)$ obtained by CSRBF methods for $n=20, 30 $ and various values of $\alpha$ are reported in Table \ref{Table .7.}.
 
 \begin{table}[h]\centering \scriptsize
\caption{CSRBFs solutions of initial slope $y'(0)$ with $N=20,30$}
\label{Table .7.}
\begin{tabular}{lcccccccc}
\hline
\multirow{2}{*}{$\alpha$} & \multicolumn{2}{c}{$Wendland_{3,5}$} & \multicolumn{2}{c}{$Wu_{3,3}$} & \multicolumn{2}{c}{$Oscillator_3$} & \multicolumn{2}{c}{$Buhman_4$}\\
\cmidrule(l){2-3} \cmidrule(l){4-5} \cmidrule(l){6-7} \cmidrule(l){8-9}
 & N=20 & N=30 & N=20 & N=30 & N=20 & N=30 & N=20 & N=30 \\
\hline
0.25 &-1.15838196 & -1.15658845 & -1.15868018 & -1.15652004 & -1.16347841 & -1.15658199 & -1.15914956 & -1.15651901 \\
0.50 & -1.19674806 & -1.19179615 & -1.19388602 & -1.19180634 & -1.19468231 & -1.19180040 & -1.19476755 & -1.19176821\\
0.75 & -1.23833094 & -1.23998794 & -1.23307537 & -1.24033353 & -1.23372756 & -1.23980995 & -1.23493307 & -1.23984689 \\
\hline
\end{tabular}
\end{table}

\begin{figure}[htbp!]
\centering \includegraphics[scale=0.4]{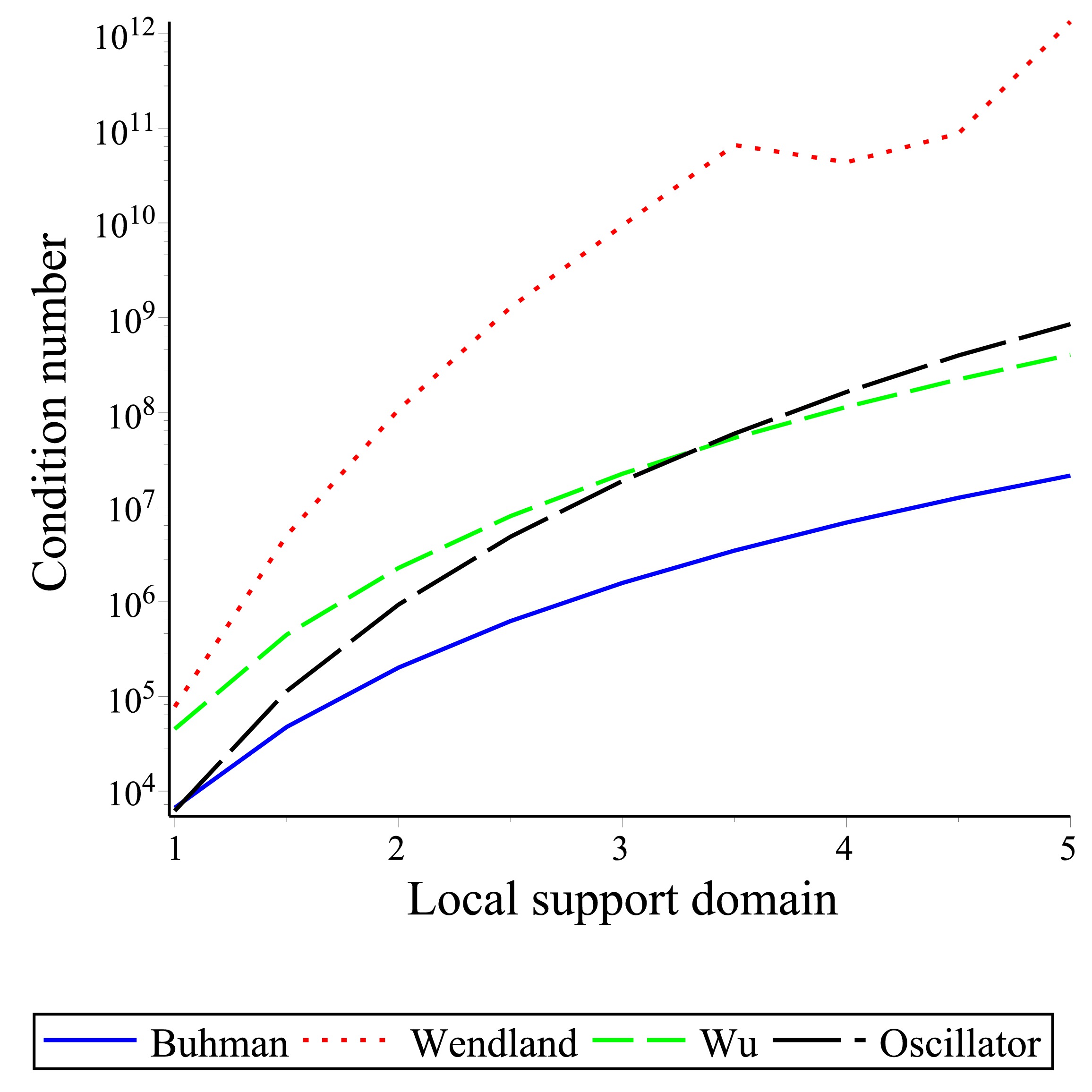}
\caption{The graph of condition number versus local support domain when $N=30$, $\alpha=0.5$ and $\varrho=1.5$}
\label{figure4}
\end{figure}

\begin{figure}[htbp!]
\centering \includegraphics[scale=0.4]{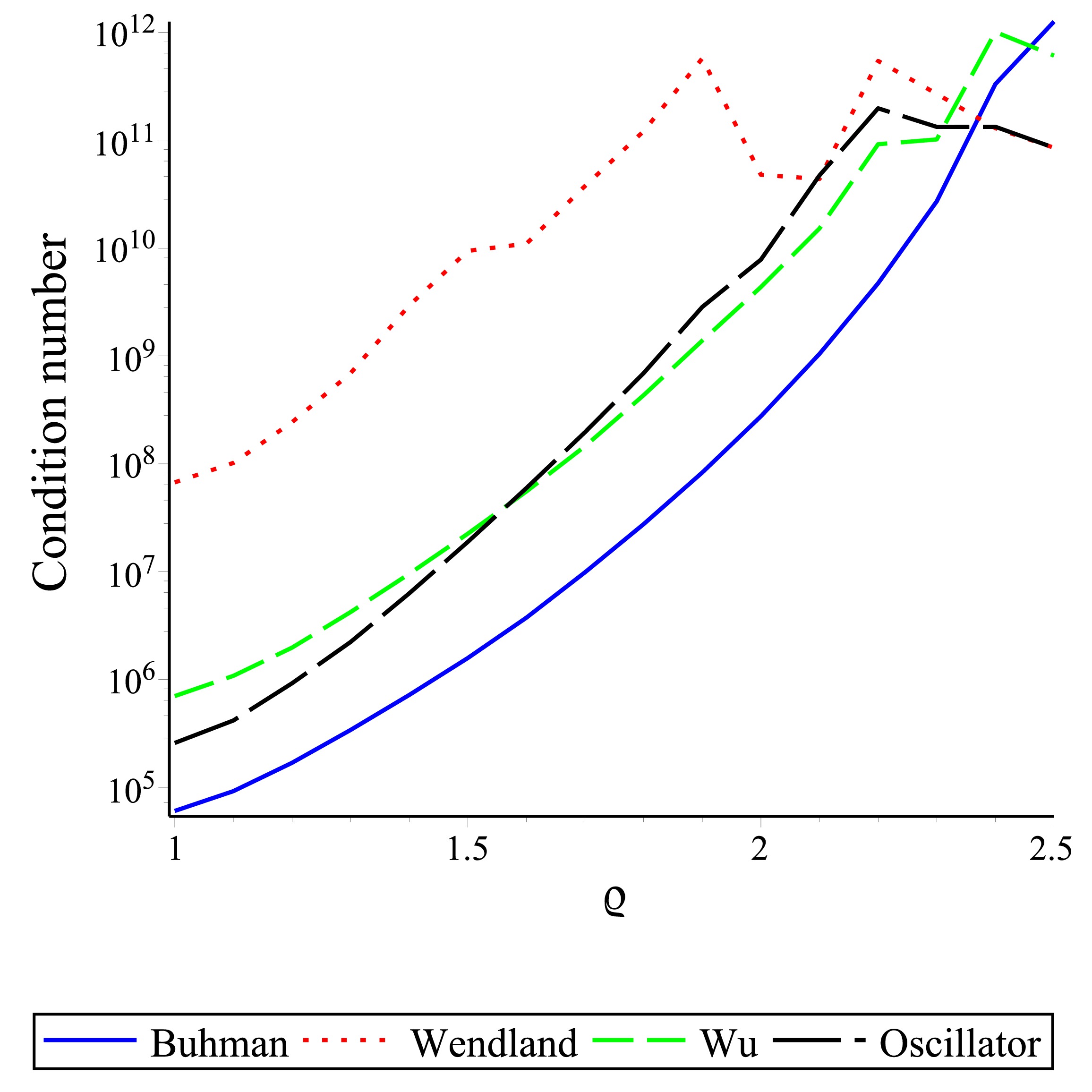}
\caption{The graph of condition number versus $\xi$ parameter when $N=30$, $\alpha=0.5$ and $r_\omega=3$}
\label{figure5}
\end{figure}

  The stability of the CSRBF scheme depends on the local support domain $r_{\omega}$. An important unsolved problem is to find a approach to determine the optimal size of $r_{\omega}$. Also, The condition number grows with $N$ and $\varrho$ for fixed values of local support domain. Figure \ref{figure4}. and \ref{figure5}. shows the condition number of matrix $A$ versus local support domain $r_{\omega}$ and $\xi$.

 Figure \ref{figure4}. shows that by decreasing the $r_\omega$, an increase is seen in the condition number of matrix $A$ and the method become more unstable. In general, the smaller the value of $r_{\omega}$, the higher percentage of zero entries in matrix $A$, so condition number matrix $A$ is a little and with increases the $r_{\omega}$, decrease percentage of zero entries in matrix $A$  and matrix $A$ will be unstable. A percentage of Zeros and condition number of CSRBFs in Matrix $A_{30\times 30}$  based on different values of local support domain $r_\omega$ illustrated in Table \ref{Table. 8.}. 
 
 \begin{table}[htbp] 
\caption{A percentage of Zero and condition number in Matrix $A_{30\times 30}$.}
\label{Table. 8.} 
\centering \scriptsize \begin{tabular}{llllll}
\hline\noalign{\smallskip}
 $r_\omega$ & $Zero$ & $Wendland_{3,5}$ & $Wu_{3,3}$ & $Buhman_4$ & $Oscillator_3$ \\
\noalign{\smallskip}\hline\noalign{\smallskip}
1.0 & 64.44\% & 7.789055e04 & 4.532194e04 & 6.667960e03 & 6.196000e03\\
1.5 & 50.45\% & 5.005580e06 & 4.490489e05 & 4.762594e04 & 1.137094e05\\
2.0 & 38.55\% & 1.070687e08 & 2.277201e06 & 2.017929e05 & 9.349735e05\\
2.5 & 28.00\% & 1.284523e09 & 8.029373e06 & 6.255295e05 & 4.865744e06\\
3.0 & 19.33\% & 5.416534e09 & 2.248322e07 & 1.581147e06 & 1.886481e07\\
3.5 & 12.00\% & 6.652775e10 & 5.373891e07 & 3.469944e06 & 5.974023e07\\
4.0 & 06.33\% & 4.377016e10 & 1.142168e08 & 6.860756e06 & 1.642163e08\\
4.5 & 02.11\% & 8.799029e10 & 2.230297e08 & 1.252852e07 & 3.990271e08\\
5.0 & 00.00\% & 1.338415e12 & 4.055797e08 & 2.148696e07 & 8.539759e08\\
\noalign{\smallskip}\hline
\end{tabular}
\end{table}

\begin{figure}[htbp!]
\centering \includegraphics[scale=0.4]{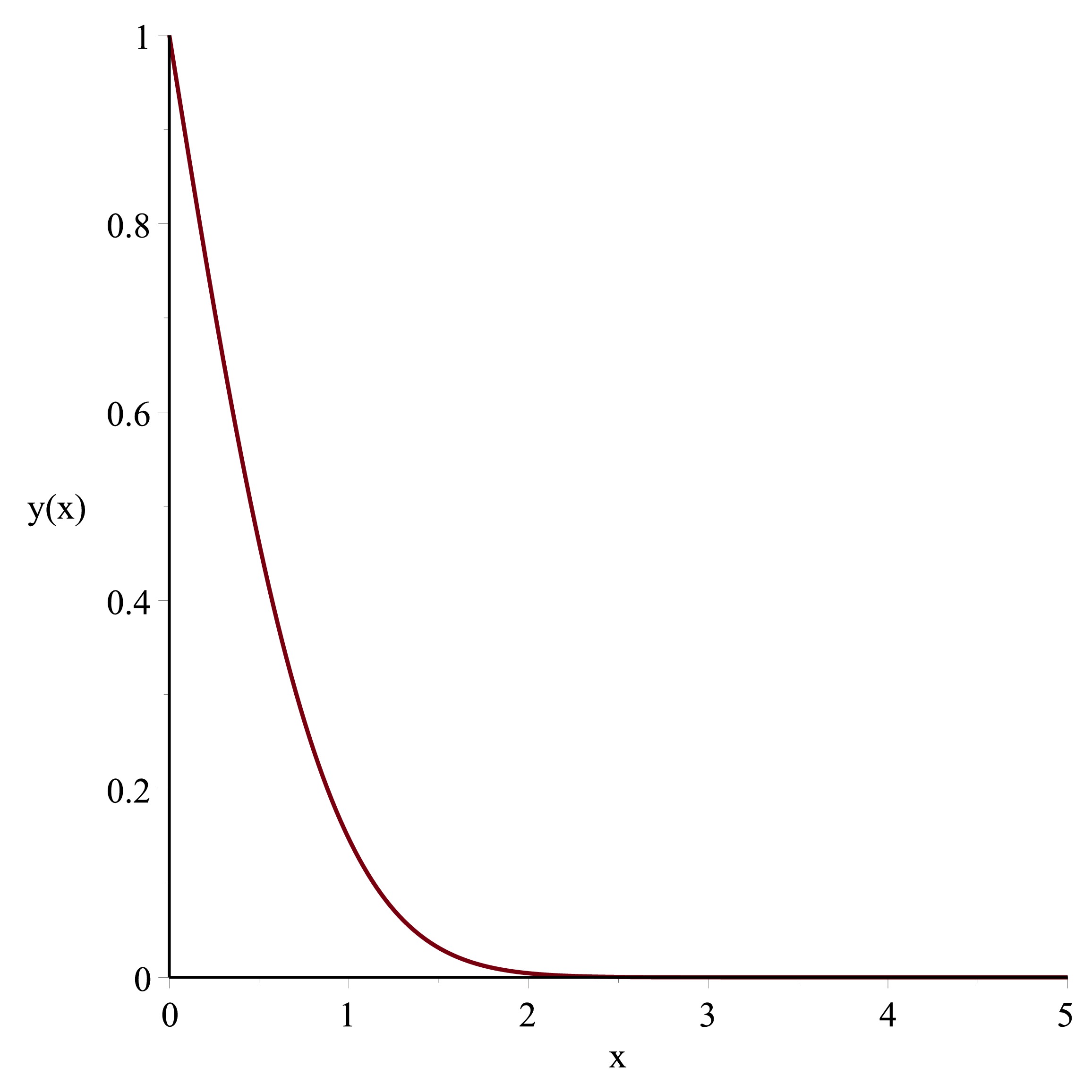}
\caption{solution of $y(x)$ when $N=30$.}
\label{figure6}
\end{figure}

 To find the critical values of $r_{\omega}$ , we use the residual error which can see in Figure \ref{figure1}.  Figure \ref{figure6}. displays the $y(x)$ for $\alpha=0.5$ in Unsteady gas equation. 
 
  \begin{figure}[htbp!]
\centering \includegraphics[scale=0.4]{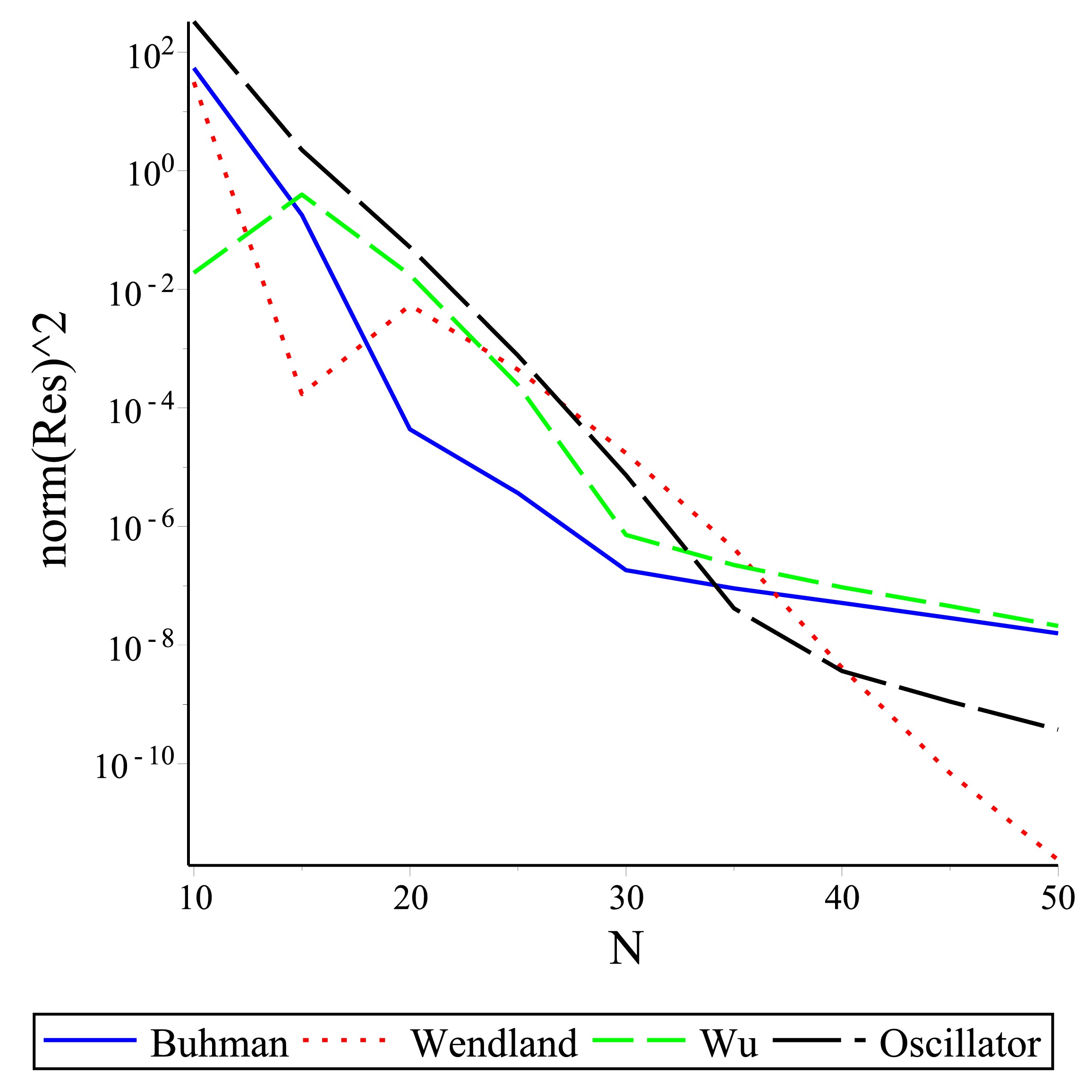}
\caption{The graph of residual error versus $N$ when $r_\omega=3$, $\alpha=0.5$ and $xi=1.5$}
\label{figure7}
\end{figure}
 
 in Figure \ref{figure7}. displays the convergence rates of the CSRBFs method for Eq. (\ref{Main_EQ}) for some value of $N$ for $\alpha=0.5$. These figures illustrate the convergence rate of the method.

\section{ Conclusion}
The fundamental goal of this paper was the construct an approximation to the solution of unsteady gas equation. a set of CSRBFs with these properties were proposed for providing an effective but simple way to improve the convergence rate. A comparison was made among the solutions of \citep{Kazem,Upad} and this work. the absolute error $\|Res\|^2$ were obtained. This paper has provided an acceptable method for Unsteady Gas equation. It was also confirmed by logarithmic figures of residual function that this method has an exponential convergence rate. Additionally, high convergence rate and good accuracy are obtained by the proposed method using relatively low numbers of collocate points.
\label{sec:6}










\end{document}